\documentclass[a4paper]{amsart}
\usepackage{amsmath,amssymb,amsthm,epic,pifont}

\newtheorem*{proposition*}{proposition}
\newtheorem{proposition}{proposition}

\newtheorem{lemma}[proposition]{Lemma}

\newtheorem{theorem}[proposition]{Theorem}

\def\AND{\qquad\text{and}\qquad}

\def\res{\operatorname{Res}}
\def\IR{\mathbb R}\def\ZZ{\mathbb Z}
\def\cint{\int_{\,\text{NW}}^{\,\text{SE}}\!d\zeta\,\,}
\def\lint{\int_{\,\text{NW}}^{\,\text{SE}}\!d\lambda\,\,}
\def\rint{\int_{-\infty}^{\,\infty}\!d\zeta\,\,}
\def\e{\mskip1mu e\mskip2mu}
\def\q{\mskip1mu q\mskip2mu}
\def\g{\mskip2mu\gamma\mskip1mu}
\def\gh{\mskip2mu\hat{\gamma\mskip1mu}}
\def\oot{{\,1/\tau}}
\def\ev{\lambda}
\def\ef{\mskip2mu\psi^{}_\ev}
\def\cf{\mskip2mu\psi^\ev}
\def\G{\mskip2mu G\mskip1mu}
\def\gddz{\g\Bigl(-\,\frac\tau{2\ipi}\,\frac d{dz}\Bigr)}

\def\x{\mskip2mu X\mskip1mu}

\def\f{\mskip2mu f\mskip1mu}
\def\om{\omega}\def\op{\omega'}
\def\opp{\omega''}
\def\ipi{i\pi}

\def\P{p}\def\Q{q}

\def\pole{\makebox(0,0){$\scriptscriptstyle\times$}}

\def\zero{\circle{3}}

\begin{document}
\title{Noncommutative Hypergeometry}
\author[A. Yu. Volkov]{Alexandre Yu. Volkov}
\address{Dienst Theoretische Natuurkunde\newline
\indent Vrije Universiteit Brussel\newline
\indent Pleinlaan 2, B-1050 Brussels, Belgium}
%\thanks{FWO}
%%%%%%%%%%%%%%%%%%%%%%%%%%%%%%%%%%%%%%%%%%%%%%%%%%%%%
\begin{abstract}
A certain special function of the generalized
hypergeometric variety is shown to fulfill a host of
useful noncommutative identities.
\end{abstract}
%%%%%%%%%%%%%%%%%%%%%%%%%%%%%%%%%%%%%%%%%%%%%%%%%%%%%
\maketitle
%%%%%%%%%%%%%%%%%%%%%%%%%%%%%%%%%%%%%%%%%%%%%%%%%%%%%
\section*{\protect\large\bf Introduction}
%%%%%%%%%%%%%%%%%%%%%%%%%%%%%%%%%%%%%%%%%%%%%%%%%%%%%
Fix complex~$\tau$ with $\operatorname{Im}\tau>0$
-- so that $\q=\e^{\ipi\tau}$ and
$\q^{-1/\tau^2}=\e^{-\ipi/\tau}$ both be less
than~$1$ in modulus --  and consider the function
\begin{equation}\label{ip}
   \g(z)=\frac{(\q^2\e^{-2\ipi z};\q^2)_\infty}
   {(\e^{-2\ipi z/\tau};\q^{-2/\tau^2})_\infty}\,,
\end{equation}
where $(a;b)_\infty$ is the usual Pochhammer-style
symbol for $(1-a)(1-ab)(1-ab^2)\ldots$,
or explicitly
$$ \g(z)=\frac{(1-\e^{-2\ipi(z-\tau)})
   (1-\e^{-2\ipi(z-2\tau)})
   (1-\e^{-2\ipi(z-3\tau)})\ldots}
   {(1-\e^{-2\ipi z/\tau})
   (1-\e^{-2\ipi(z+1)/\tau})
   (1-\e^{-2\ipi(z+2)/\tau})\ldots}\,,             $$
Clearly, this  is meromorphic at all
$z\neq\infty$ and has a remarkable pattern of zeros
and poles:
%%%%%%%%%%%%%%%%%%%%%%%%%%%%%%%%%%%%%%%%%%%%%%%%%%%%%
\begin{figure}[h]\centering
%\framebox{
 \begin{picture}(300,92)
  \put(50,28){\begin{picture}(0,0)
   \put(14,-28){\line(-1,2){34}}
   \put(-45,0){\line(1,0){75}}
%   \put(40,-10){\makebox(0,0)[br]{$\IR$}}
   \put(-6,12){\circle*{3}}
   \put(0,13){\makebox(0,0)[l]{$\tau$}}
   \put(10,0){\circle*{3}}
   \put(16,-4){\makebox(0,0)[t]{$1$}}
   \multiput(-18,36)(-10,0){3}{\zero}
   \multiput(-12,24)(-10,0){4}{\zero}
   \multiput(-6,12)(-10,0){4}{\zero}
   \multiput(0,0)(-10,0){5}{\zero}
   \multiput(6,-12)(-10,0){5}{\zero}
   \multiput(12,-24)(-10,0){6}{\zero}
  \end{picture}}
  \put(100,52){\begin{picture}(0,0)
   \put(14,-28){\line(-1,2){34}}
   \put(-35,0){\line(1,0){75}}
%   \put(40,-10){\makebox(0,0)[br]{$\IR$}}
   \put(-6,12){\circle*{3}}
   \put(-8,6){\makebox(0,0)[r]{$\tau$}}
   \put(10,0){\circle*{3}}
   \put(16,-4){\makebox(0,0)[t]{$1$}}
   \multiput(-18,36)(-10,0){2}{\zero}
   \multiput(-12,24)(-10,0){3}{\zero}
   \multiput(-6,12)(-10,0){3}{\zero}
   \multiput(-8,36)(10,0){5}{\zero}
   \multiput(-2,24)(10,0){4}{\zero}
   \multiput(4,12)(10,0){4}{\zero}
  \end{picture}}
  \put(235,40){\begin{picture}(0,0)
   \put(20,-40){\line(-1,2){46}}
   \put(-55,0){\line(1,0){115}}
   \put(-34,48){\makebox(0,0){N}}
   \put(28,-36){\makebox(0,0){S}}
   \put(55,-4){\makebox(0,0)[t]{E}}
   \put(-50,5){\makebox(0,0)[b]{W}}
%   \put(60,-10){\makebox(0,0)[br]{$\IR$}}
   \put(-6,12){\circle*{3}}
   \put(-11,13){\makebox(0,0)[r]{$\tau$}}
   \put(10,0){\circle*{3}}
   \put(16,-4){\makebox(0,0)[t]{$1$}}
   \multiput(-14,48)(10,0){7}{\zero}
   \multiput(-8,36)(10,0){7}{\zero}
   \multiput(-2,24)(10,0){6}{\zero}
   \multiput(4,12)(10,0){6}{\zero}
   \multiput(0,0)(-10,0){6}{\pole}
   \multiput(6,-12)(-10,0){6}{\pole}
   \multiput(12,-24)(-10,0){7}{\pole}
   \multiput(18,-36)(-10,0){7}{\pole}
%   \qbezier[22](-53,6)(-10,0)(-3,6)
%   \qbezier[15](-3,6)(4,12)(5,0)
%   \qbezier[21](5,0)(6,-12)(26,-36)
   \qbezier[15](-47,36)(-35,30)(-29,18)
   \qbezier[20](-29,18)(-23,6)(-3,6)
   \qbezier[20](-3,6)(17,6)(23,-6)
   \qbezier[15](23,-6)(29,-18)(41,-24)
   \put(43,-25){\vector(2,-1){0}}
%%%   \qbezier[19](-54,8)(-13,6)(-3,6)
%%%   \qbezier[6](-3,6)(2,6)(5,0)
%%%   \qbezier[18](5,0)(11,-12)(26,-36)
%%%   \put(-33,7){\vector(-1,0){0}}
%   \qbezier[32](-55,30)(5,30)(38,-36)
%%   \qbezier[27](-54,18)(6,18)(33,-36)
%%   \put(-31,16){\vector(-4,1){0}}
%   \put(-56,30){\vector(-1,0){0}}
%   \put(-19,23){\vector(-3,1){0}}
%   \put(18,-28){\vector(2,-3){0}}
%   \put(-18,3){\vector(-1,0){0}}
%   \put(10,-13){\vector(-1,2){0}}
%   \dottedline{2}(-50,19)(46,-5)
%   \put(46,-5){\vector(4,-1){0}}
 \end{picture}}\end{picture}
%}
\caption{The denominator, numerator and whole of
function~$\g$. As an exercise, figure out what
happens as~$\tau$ approaches the real line. Will it
matter on which half it lands?}
\end{figure}
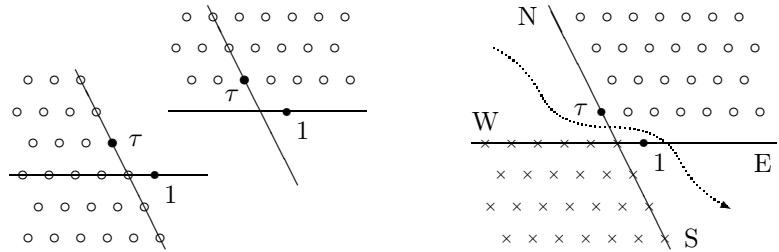
%%%%%%%%%%%%%%%%%%%%%%%%%%%%%%%%%%%%%%%%%%%%%%%%%%%%%
all are simple, and, as shown on the figure, they
fill the north-eastern and south-western quarters --
$k,l>0$ and $k,l\leq0$ respectively --
of the lattice~$k+l\tau$. With slight variations,
this function, called below the $\g$-function, has of
late been circulating in connection with quantum
integrable models under fancy names like double sine
or noncompact quantum dilogarithm.\footnote
{See (Faddeev Kashaev Volkov 2001),
(Kharchev Lebedev Semenov-Tian-Shansky 2003),
(Bytsko Teschner 2003) and references therein.
Ultimately, though, this function traces
back to (Shintani 1977) and (Barnes 1899).}
Yet it deserves more attention than has so far been
given, and so this article aims to introduce the
matter to a wider mathphysical audience, and tie up
a few loose ends in the process.
%%%%%%%%%%%%%%%%%%%%%%%%%%%%%%%%%%%%%%%%%%%%%%%%%%%%%
\section{\protect\large\bf Two equations}
%%%%%%%%%%%%%%%%%%%%%%%%%%%%%%%%%%%%%%%%%%%%%%%%%%%%%
Note that our function satisfies difference equations
\begin{equation}\label{te}
   \frac{\g(z+\tau)}{\g(z)}=1-\e^{-2\ipi z}\AND
   \frac{\g(z+1)}{\g(z)}=1-\e^{-2\ipi z/\tau},
\end{equation}
and is indeed their only common solution -- up to
multiplication by an elliptic function of periods~$1$
and~$\tau$ that is, but this little ambiguity is
easily avoided by adding a kind of minimal asymptotic
condition that $\g(z)$ goes to~$1$ as $z$ goes
southeast, that is $z\rightarrow\infty$ in the sector
$\arg(-\tau)<\arg z<0$. So, these two equations can
and will serve as a workable alternative to the
infinite product definition. They also help explain
the rather obscure title of this article.

On the one hand, each of the two resembles
the most famous difference equation
$$ \frac{\Gamma(z+1)}{\Gamma(z)}=z\,,              $$
which defines Euler's Gamma function. Hence
`hypergeometry', and hence our function ought to be
called, say, a double or elliptic gamma function --
but since those names are already taken, let us
settle for $\g$-function.

On the other hand, our two equations can be rewritten
as
$$ \e^{-2\ipi z}\g(z)+\g(z+\tau)=\g(z)
   \AND\e^{-2\ipi z/\tau}\g(z)+\g(z+1)=\g(z)\,,    $$
or in operator form
\begin{equation*}\tag{\ref{te}$'$}
   (\e^{-2\ipi z}+\e^{\tau\,d/dz})\g=\g\AND
   (\e^{-2\ipi z/\tau}+\e^{d/dz})\g=\g\,,
\end{equation*}
where, by usual abuse of notation,~$z$ and
$d/dz$ stand for operators of multiplication and
differentiation by~$z$, that is
$$ (zf)(z)=zf(z)\AND\frac d{dz}f(z)=f'(z)\,.       $$
This identifies the $\g$-function as an eigenfunction
of useful operators and so gives everything we do
some sort of `noncommutative' meaning. For instance,
if those operators share one eigenfunction, might
they also share the rest of them and be therefore
functions of each other? Yes indeed, it is easily
checked that functions
$$ \ef(z)=\e^{-2\ipi\ev z/\tau}\g(z-\ev)           $$
satisfy full spectral equations
\begin{equation}\label{sp}
   (\e^{-2\ipi z}+\e^{\tau\,d/dz})\ef
   =\e^{-2\ipi\ev}\ef\quad\text{and}\quad
   (\e^{-2\ipi z/\tau}+\e^{d/dz})\ef
   =\e^{-2\ipi\ev/\tau}\ef\,.
\end{equation}
Hence, by comparing the respective eigenvalues and
optimistically assuming that functions $\ef$ span
some reasonable functional space like~$L^2$ on the
dotted line in Figure~1, follows a somewhat
surprising relation
\begin{equation}\label{oot}
   (\e^{-2\ipi z}+\e^{\tau\,d/dz})^\oot
   =\e^{-2\ipi z/\tau}+\e^{d/dz},
\end{equation}
which has previously only been noticed in one
particular case that $1/\tau$ is positive integer.
Let us not get ahead of ourselves, though, and get
over with the hypergeometric part first.
%%%%%%%%%%%%%%%%%%%%%%%%%%%%%%%%%%%%%%%%%%%%%%%%%%%%%
\section{\protect\large\bf Reflection formula}
%%%%%%%%%%%%%%%%%%%%%%%%%%%%%%%%%%%%%%%%%%%%%%%%%%%%%
Recall the classical formula
$$ \Gamma(z)\Gamma(1-z)=\frac{\pi}{\sin\pi z}\,,   $$
which says that reflection about the point~$1/2$
reduces the Gamma function to an elementary one. The
same happens to be true of the $\g$-function except
the natural reflection point is now $(1+\tau)/2$.
Indeed, since poles and zeros of the $\g$-function
are symmetric to each other about the said point,
the product
$$ \G(z)=\g(z)\g(1+\tau-z)                         $$
has none of either. It might then equal a Gaussian
exponential, and so it turns out. By one of
equations~(\ref{te}) we have
$$ \frac{\G(z+\tau)}{\G(z)}
   =\frac{\g(z+\tau)\g(1-z)}{\g(z)\g(1+\tau-z)}
   =\frac{1-\e^{-2\ipi z}}{1-\e^{-2\ipi(\tau-z)}}
   =-\e^{-2\ipi z},                                $$
and by the other
$\G(z+1)/\G(z)=\cdots=-\e^{2\ipi z/\tau}$.
Since the same equations are solved by
$\e^{\ipi z(1+\tau-z)/\tau}$, we have
$$ \G(z)=\e^{\ipi z(1+\tau-z)/\tau}
   \times\text{an elliptic function with
   periods $1$ and $\tau$}\,,                      $$
but without zeros nor poles on either side, that
elliptic factor can only be constant. Set $z=1$ to
show that that constant equals $-\g(1)\g(\tau)$, then
set $z=0$ in equations~(\ref{te}) to show that
$$ \g(\tau)=2\ipi\res\g(z)|_{z=0}=\tau\g(1)\,.     $$
Hence the anticipated `reflection formula':
\begin{equation}\label{rf}
   \g(z)\g(1+\tau-z)=-\,\tau\g(1)^2\,
   \e^{\ipi z(1+\tau-z)/\tau}.
\end{equation}
This is good enough in this context, yet the question
remains whether $\g(1)$ could be evaluated in
absolute terms. Those familiar with Dedekind's eta
function should already know because
$$ \g(1)=\frac{(\q^2;\q^2)_\infty}
   {(\q^{-2/\tau^2};\q^{-2/\tau^2})_\infty}
   =\frac{\e^{-\ipi\tau/12}\,\eta(\tau)}
   {\e^{\ipi/12\tau}\,\eta(-1/\tau)}\,,            $$
but the rest of us will have to wait until Section~4.
%%%%%%%%%%%%%%%%%%%%%%%%%%%%%%%%%%%%%%%%%%%%%%%%%%%%%
\section{\protect\large\bf Under Fourier transform}
%%%%%%%%%%%%%%%%%%%%%%%%%%%%%%%%%%%%%%%%%%%%%%%%%%%%%
Recall another classical formula called Euler's Gamma
integral, which reads
$$ \Gamma(z)=\int_0^\infty
   \frac{dw}w\,\,w^{\,z}\e^{-w},                   $$
and says that Fourier transform (or Mellin to be
precise) reduces the Gamma function to the
exponential one. The same cannot be quite true of the
$\g$-function, for, as we remember, its defining
equations~(\ref{te}$'$) mix differentiation and
multiplication operators in a symmetric manner. So,
since Fourier transform maps those operators more or
less into each other, it will map the
$\g$-function more or less into itself rather than
reduce it to anything else. Specifically, consider
the Fourier integral
$$ \gh(z)=-\,\frac1\tau\cint
   \e^{-2\ipi \zeta z/\tau}\g(\zeta)               $$
along that same dotted line in Figure~1. Now one of
equations~(\ref{te}$'$) translates into
$$ \gh(z)=\frac1\tau\cint\e^{-2\ipi\zeta z/\tau}
   (\e^{-2\ipi\zeta/\tau}+\e^{d/d\zeta})\g(\zeta)
   =(\e^{d/dz}+\e^{2\ipi z/\tau})\gh(z)\,,         $$
and the other into
$\gh=\cdots=(\e^{\tau\,d/dz}+\e^{2\ipi z})\gh(z)$,
which are, as expected, virtually the same equations
as~(\ref{te}$'$) themselves. Clearly, their general
solution is
$$ \gh(z)=\frac{\text{an elliptic function with
periods $1$ and $\tau$}}{\g(1+\tau-z)}\,,          $$
but since our integral converges too well for
elliptic functions to creep in, that elliptic factor
is again constant -- and equals, of course,
$(2\ipi/\tau)\*\res\g(z)|_{z=0}$, that is, as we
already know, that same $\g(1)$. Hence, after
inversion, follows the `tau-gamma integral'
\begin{equation}\label{ft}
   \g(z)=\g(1)\cint\frac
   {\e^{2\ipi\zeta z/\tau}}{\g(1+\tau-\zeta)}\,,
\end{equation}
which confirms that the $\g$-function is indeed a
Fourier image of more or less itself rather than some
quasiexponential function. Or is it that the
$\g$-function somehow emulates both the gamma and
exponential functions? We will find out in Section~6.
%%%%%%%%%%%%%%%%%%%%%%%%%%%%%%%%%%%%%%%%%%%%%%%%%%%%%
\section{\protect\large\bf Tau-binomial theorem}
%%%%%%%%%%%%%%%%%%%%%%%%%%%%%%%%%%%%%%%%%%%%%%%%%%%%%
The more immediate question is, what about the beta
integral? This is settled by
Kashaev-Ponsot-Teschner's `tau-binomial theorem':
for all~$y\neq k+l\tau$ ($\ZZ\ni k,l\leq0$) we have
\begin{equation}\label{tbt}
   \frac{\g(y)\g(z)}{\g(y+z)}
   =\g(1)\cint\frac{\e^{2\ipi\zeta z/\tau}
   \g(y-\zeta)}{\g(1+\tau-\zeta)}\,,
\end{equation}
provided the middle part of the integration line is
rerouted, if necessary, so as to separate zeros of
$\g(1+\tau-\zeta)$ from poles of
$\g(y-\zeta)$.\footnote
{The same rule applies, without further mention, to
all the integrals below: the integration line must
always separate the southwestern and northeastern
`sequences' of poles.}
This is a straightforward extension of the tau-gamma
integral~(\ref{ft}), in the sense that it reduces
to the latter as~$y$ goes southeast, and is verified
along the same lines, that is by comparing equations
with respect to~$z$ and evaluating contribution of
the pole at $\zeta=0$. The next question, then, is
whether further extension is possible. The short
answer is no, there are no more Fourier integrals
left to take. A longer answer follows in the next
section, but to wrap this one up, let us consider
another two useful limit cases.

First, send $y+z$ to zero to obtain in the limit
$$ \delta(y)=-\,\frac1\tau\cint
   \frac{\e^{-2\ipi\zeta y/\tau}
   \g(y-\zeta)}{\g(1+\tau-\zeta)}\,,              $$
or, after a suitable change of variables,
\begin{equation}\label{delta}
   \delta(z-y)=\lint\cf(z)\ef(y)\,,
\end{equation}
where $\ef(y)=\e^{-2\ipi\zeta y/\tau}\g(y-\zeta)$
are those same (generalized) eigenfunctions from
Section~1, and
$\cf(z)=(-1/\tau)\,\e^{2\ipi\ev z/\tau}/
\g(1+\tau-\ev+z)$.
I leave it to the reader to figure out the details,
but the upshot is, anyway, that functions~$\ef$ do
indeed span a wide range of functional spaces, which
consist of functions that are, loosely speaking,
well-behaved in the northwestern and southeastern
quarters (as mapped in Figure~1) of the complex
plane. This bodes well for potential `noncommutative'
applications, but, again, let us get over with the
hypergeometric part first.

Second, apply the reflection formula~(\ref{rf}) a few
times and change variables so that the tau-binomial
theorem~(\ref{tbt}) becomes
\begin{equation*}\tag{\ref{tbt}$'$}
   \frac{\g(1+\tau-z)\g(y)}{\g(1+\tau-z+y)}
   =-\,\tau\g(1)^3\cint
   \frac{\e^{\ipi\zeta(1+\tau-\zeta)/\tau}}
   {\g(1+\tau-z+\zeta)\g(1+\tau-\zeta+y)}\,,
\end{equation*}
then send both~$y$ and~$-z$ southeast to reduce
things to the Gaussian integral
$$ 1=-\,\tau\g(1)^3\cint
   \e^{\ipi\zeta(1+\tau-\zeta)/\tau}\,.            $$
Hence an opportunity to evaluate~$\g(1)$ without
Dedekind's help. Take the integral and then, in order
to pick the right cubic root, go back to the infinite
product expansion~(\ref{ip}) and see that $\g(1)$
should equal~$1$ if $\tau=i$. Thus,
\begin{equation}\label{ded}
   \g(1)=\,i^{\,5/6}\,\tau^{\,-1/2}\,
   \e^{-\ipi(1+\tau)^2/12\tau},
\end{equation}
but I am not sure what to make of it, and will be
using~$\g(1)$ for shorthand anyway.
%%%%%%%%%%%%%%%%%%%%%%%%%%%%%%%%%%%%%%%%%%%%%%%%%%%%%
\section{\protect\large\bf Beyond Fourier transform}
%%%%%%%%%%%%%%%%%%%%%%%%%%%%%%%%%%%%%%%%%%%%%%%%%%%%%
Write
$$ \frac{\g(z-\zeta)}{\g(1+\tau-\zeta)}
   =\frac{\g(y+z-\zeta)}{\g(1+\tau-\zeta)}\,
   \frac{\g(z-\zeta)}{\g(y+z-\zeta)}               $$
and take Fourier transform of both sides (using the
tau-binomial theorem~(\ref{tbt}) twice directly and
once in reverse, and the fact that Fourier transform
of a product is a convolution of Fourier transforms
of its factors) to obtain
$$ \frac{\g(x)\g(y)\g(z)}{\g(x+z)\g(y+z)}
   =\g(1)\cint\frac{\e^{2\ipi\zeta z/\tau}
   \g(x-\zeta)\g(y-\zeta)}
   {\g(x+y+z-\zeta)\g(1+\tau-\zeta)}\,.            $$
This looks very like the $\tau$-binomial
theorem~(\ref{tbt}) and obviously reduces to the
latter as either~$x$ or $y$ go southeast -- yet the
right hand side is no longer a Fourier integral of
course. So, the longer answer to the question of the
previous section is that extension of the
tau-binomial theorem~(\ref{tbt}) is possible after
all, but it turns out to be an `addition theorem'
rather than Fourier integral. How about another
extension then?

Use the reflection formula~(\ref{rf}) a few times and
change variables so that the above formula takes a
more transparent form
\begin{multline}\label{gauss}
   \frac{\g(\nu+\mu)\g(\nu+\kappa)\g(\mu)
   \g(\kappa)}{\g(\nu+\mu+\kappa)}                 \\
   =-\,\frac1{\tau\g(1)}\cint
   \e^{2\ipi\zeta(\nu+\zeta)/\tau}\g(\nu+\zeta)
   \g(\mu-\zeta)\g(\kappa-\zeta)\g(\zeta)\,,
\end{multline}
in which it unmistakably resembles the addition
theorem for binomial coefficients\footnote
{Just in case, this expresses equality of
coefficients of $w^m$ in both hand sides of the
formula $(1+w)^{n+m+k}=(1+w)^{n+m}(1+w)^k$.}
$$ \frac{(n+m+k)!}{(n+m)!(n+k)!m!k!}
   =\sum_j\frac1{(n+j)!(m-j)!(k-j)!j!}\,.          $$
Then, since the latter is known to have exactly one
extension in the shape of Pfaff-Saalsch\"utz's sum
$$ \frac{(n+l+k)!(n+m+k)!(n+m+l)!}
   {(n+m)!(n+l)!(n+k)!m!l!k!}
   =\sum_j\frac{(n+m+l+k-j)!}
   {(n+j)!(m-j)!(l-j)!(k-j)!j!}\,.                 $$
it is a safe educated guess that the $\g$-function
satisfies a similar `ultimate integral identity'
\begin{multline}\label{pfaff}
   \frac{\g(\nu+\mu)\g(\nu+\lambda)
   \g(\nu+\kappa)\g(\mu)\g(\lambda)\g(\kappa)}
   {\g(\nu+\lambda+\kappa)\g(\nu+\mu+\kappa)
   \g(\nu+\mu+\lambda)}                            \\
   =-\,\frac1{\tau\g(1)}\cint
   \e^{2\ipi\zeta(\nu+\zeta)/\tau}
   \frac{\g(\nu+\zeta)\g(\mu-\zeta)\g(\lambda-\zeta)
   \g(\kappa-\zeta)\g(\zeta)}
   {\g(\nu+\mu+\lambda+\kappa-\zeta)}
\end{multline}
-- at least safe enough not to bother verifying it
just yet. We will do it anyway in the noncommutative
part, which now begins.
%%%%%%%%%%%%%%%%%%%%%%%%%%%%%%%%%%%%%%%%%%%%%%%%%%%%%
\section{\protect\large\bf Going noncommutative}
%%%%%%%%%%%%%%%%%%%%%%%%%%%%%%%%%%%%%%%%%%%%%%%%%%%%%
Let us go back to the tau-binomial theorem to try and
interpret it as an operator relation. The
form~(\ref{tbt}$'$) is best suited for that. Apply
the reflection formula~(\ref{rf}) one more time to
obtain
\begin{equation*}\tag{\ref{tbt}$''$}
   \g(1+\tau-z)\gh(z-y)\g(y)=\cint\gh(z-\zeta)
   \g(1+\tau-\zeta)\g(\zeta)\gh(\zeta-y)\,,
\end{equation*}
where, as before,~$\gh(z)=\g(1)/\g(1+\tau-z)$ is
Fourier image of the $\g$-function. This indeed lends
itself to be interpreted as operator relation
$$ \mathit{abc=bacb}\,,                            $$
where~$a$ and~$c$ are operators of pointwise
multiplication by $\g(1+\tau-z)$ and $\g(z)$,
$$ a\f(z)=\g(1+\tau-z)\f(z)\AND c\f(z)=\g(z)\f(z), $$
and~$b$ is that of convolution with $\gh$,
$$ b\f(z)=\cint\gh(z-\zeta)\f(\zeta).              $$
So if~$z$ and $d/dz$ stand, as before, for operators
of multiplication and differentiation by~$z$, then,
by near tautology,
$$ a=\g(1+\tau-z)\AND c=\g(z)\,,                   $$
and, by a textbook argument about multiplication vs
convolution,
\begin{multline*}
   \cint\gh(z-\zeta)\f(\zeta)
   =\cint\gh(\zeta)\f(z-\zeta)                     \\
   =\cint\gh(\zeta)\e^{-\zeta\,d/dz}\f(z)
   =\gddz\f(z)\,,
\end{multline*}
that is
$$ b=\gddz\,.                                      $$
Thus, a direct operator translation of the
tau-binomial theorem reads
\begin{equation}\label{seven}
   \g(1+\tau-z)\gddz\g(z)
   =\gddz\g(1+\tau-z)\g(z)\gddz\,,
\end{equation}
but before you say a word, here is another, not
so direct but shorter one.

Apply the reflection formula~(\ref{rf}) to a
different part
of~(\ref{tbt}$'$) to obtain
\begin{equation*}\tag{\ref{tbt}$'''$}
   \gh(z-y)\g(y)=\g(z)\cint
   \e^{\ipi(\zeta-z)(1+\tau-\zeta-z)/\tau}
   \gh(z-\zeta)\gh(\zeta-y)\,,
\end{equation*}
or
$$ \mathit{bc=ceb}\,,                              $$
where operators~$b$ and~$c$ are the same as above,
and~$e$ acts as
$$ e\f(z)=\cint\e^{\ipi(\zeta-z)
   (1+\tau-\zeta-z)/\tau}\gh(z-\zeta)\f(\zeta)\,.  $$
It is then easy to figure out that
$$ e=\g\Bigl(-\,\frac\tau{2\ipi}\,\frac d{dz}
   +z-\frac{1+\tau}2\Bigr)                         $$
and thus obtain Kashaev's
`pentagon identity'\footnote
{No, it is not called that because it features five
factors. See (Kashaev 2000) for an explanation and
further references.}
\begin{equation}\label{pentagon}
   \gddz\g(z)=\g(z)\g\Bigl(-\,\frac\tau{2\ipi}\,
   \frac d{dz}+z-\frac{1+\tau}2\Bigr)\gddz\,.
\end{equation}

If this is still not good enough,
rewrite~(\ref{tbt}) as
\begin{equation*}\tag{\ref{tbt}$''''$}
   \g(\ev)\ef(z)
   =\g(z)\cint\gh(z-\zeta)\ef(\zeta)\,,
\end{equation*}
or
$$ \g(\ev)\ef=\mathit{cb}\ef\,,                    $$
where~$b$ and~$c$ are the same as above, and
$\ef(z)=\e^{-2\ipi\ev z/\tau}\g(z-\ev)$ are the
same eigenfunctions that have already appeared twice
on these pages (Sections~1 and~4). So, this is just
another spectral equation on functions that already
satisfy two. Hence, by comparing the respective
eigenvalues here and, say, in the first of
equations~(\ref{sp}), we have
$$ \g\Bigl(-\,\frac{\log(\e^{-2\ipi z}
   +\e^{\tau\,d/dz})}{2\ipi}\Bigr)=\g(z)\gddz\,,   $$
or, more compactly,
\begin{equation}\label{sb}
   \x(\e^{-2\ipi z}+\e^{\tau\,d/dz})
   =\x(\e^{-2\ipi z})\,\x(\e^{\tau\,d/dz})\,,
\end{equation}
where~$\x$ is the function such that
$$ \x(\e^{-2\ipi z})=\g(z)\,,                      $$
that is
\begin{multline*}
   \x(w)=\g\Bigl(-\,\frac{\log w}{2\ipi}\Bigr)
   =\frac{(\q^2w;\q^2)_\infty}
   {(w^\oot;\q^{-2/\tau^2})_\infty}                \\
   =\frac{(1-\q^2w)(1-\q^4w)(1-\q^6w)\ldots}
   {\bigl(1-w^\oot\bigr)
   \bigl(1-\q^{-2/\tau^2}w^\oot\bigr)
   \bigl(1-\q^{-4/\tau^2}w^\oot\bigr)\ldots}\,.
\end{multline*}
This is called `Sch\"utzenberger's equation'
after the famous French combinatorialist who
discovered such noncommutative exponentiality fifty
years ago. He did without those scary $1/\tau$-th
powers though. I will explain after a remark.

Remember we were wondering how the same function
could emulate the Gamma and exponential functions
at the same time? Now we know. The $\x\!$-function
may have zeros, poles and a cut, but it is the
exponential property that counts, and on this
grounds alone it should be accepted as a legitimate
`noncommutative' exponential function. It must be
stressed, though, that the exponential property
itself has also become `noncommutative'. What
happens, then, if the same factors are multiplied
the other way around? As it turns out, this: 
\begin{equation}\label{cc}
   \x(\e^{\tau\,d/dz})\,\x(\e^{-2\ipi z})
   =\x(\e^{-2\ipi z}-\e^{-2\ipi z}\e^{\tau\,d/dz}
   +\e^{\tau\,d/dz})\,.
\end{equation}
The derivation is quite similar to that of
Sch\"utzenberger's equation~(\ref{sb}) and is
therefore left as an exercise.
%%%%%%%%%%%%%%%%%%%%%%%%%%%%%%%%%%%%%%%%%%%%%%%%%%%%%
\section{\protect\large\bf Breakdown}
%%%%%%%%%%%%%%%%%%%%%%%%%%%%%%%%%%%%%%%%%%%%%%%%%%%%%
An explanation is indeed in order and not just of
where Sch\"utzenberger fits in all this, but, more
broadly, of how that all-important exponential
property~(\ref{sb}) could be derived step by step
rather than, as above, pulled out of the hat. Let us,
then, rederive it in that sort of heuristic
$q$-algebraic style typical of the subject.

Recall (or rederive) that operators $\e^{-2\ipi z}$
and $\e^{\tau\,d/dz}$ satisfy Weyl's relation
$$ \e^{-2\ipi z}\e^{\tau\,d/dz}
   =\q^2\,\e^{\tau\,d/dz}\e^{-2\ipi z}             $$
and, for now, forget all else. That is, consider
instead formal operators~$u$ and~$v$ only subject to
the relation
\begin{equation*}\tag{i}
   uv=\q^2vu\,.
\end{equation*}
Apply the latter repeatedly to show that
\begin{multline*}
   (u-uv+v)(1-\q^2v)(1-\q^4v)(1-\q^6v)\ldots       \\
   =(1-\q^2v)(u-\q^2uv+v)(1-\q^4v)(1-\q^6v)\ldots  \\
   =(1-\q^2v)(1-\q^4v)(u-\q^4uv+v)(1-\q^6v)\ldots  \\
   =\cdots=(1-\q^2v)(1-\q^4v)(1-\q^6v)\ldots(u+v)\,,
\end{multline*}
that is
$(u-uv+v)\*(\q^2v;\q^2)_\infty
=(\q^2v;\q^2)_\infty\*(u+v)$,
and similarly that
$(u+v)\*(\q^2u;\q^2)_\infty
=(\q^2u;\q^2)_\infty\*(u-uv+v)$,
and therefore
\begin{multline*}
   (u+v)(\q^2u;\q^2)_\infty(\q^2v;\q^2)_\infty
   =(\q^2u;\q^2)_\infty(u-uv+v)(\q^2v;\q^2)_\infty \\
   =(\q^2u;\q^2)_\infty(\q^2v;\q^2)_\infty(u+v)\,.
\end{multline*}
Thus, $(\q^2u;\q^2)_\infty(\q^2v;\q^2)_\infty$
commutes with $u+v$, and must therefore be its
function. Call it~$F$ and set $u=0$ and $v=w$ --
which is permitted by relation~(i) of course --
to see that
$$ F(w)=F(0+w)=(0;\q^2)_\infty(\q^2w;\q^2)_\infty
   =(\q^2w;\q^2)_\infty\,.                         $$
Hence what actually was Sch\"utzenberger's
discovery:\footnote
{By the way, the crucial group-likeness property of
Drinfel'd's universal $sl_2$ R-matrix,
$$ \Delta\otimes\operatorname{id}(R)=R^{13}R^{23}\AND
   \operatorname{id}\otimes\Delta(R)=R^{13}R^{12}, $$
is really just this formula in fancy disguise. See
(Faddeev 2000) and (Bytsko Teschner 2003) for more on
the Quantum Group connection.}
\begin{equation*}%\label{rsb}
   (\q^2(u+v);\q^2)_\infty
   =(\q^2u;\q^2)_\infty(\q^2v;\q^2)_\infty\,,
\end{equation*}
or in words, the numerators alone already satisfy
Sch\"utzenberger's equation.

Turning to the denominators, we seem to be stuck,
because formal operators may only be raised to
positive integer powers -- which $1/\tau$ is not.
Still, note that for positive integers we have
$$ u^mv^n=\q^{2mn}v^nu^m,                          $$
and assume, for lack of a better idea, that this
somehow remains true if one or both powers are no
longer integer. Then, whatever $u^\oot$ and $v^\oot$
might really be, they are bound, on one hand,
to satisfy Weyl's relation with $\q^{-1/\tau^2}$
instead of~$\q$,
\begin{equation*}\tag{ii}
   v^\oot u^\oot
   =\q^{-2/\tau^2}u^\oot v^\oot,
\end{equation*}
and on the other, to commute with~$u$ and~$v$:
\begin{equation*}\tag{iii}
   uv^\oot=\q^{2/\tau}v^\oot u
   =\e^{2\ipi\tau/\tau}v^\oot u=v^\oot u
   \quad\text{and}\quad
   u^\oot v=\cdots=vu^\oot.
\end{equation*}
But what about $(u+v)^\oot$ then? Note that
$$ (uv^{-1})(u+v)=\q^2(u+v)(uv^{-1}),              $$
and therefore, by the same little trick that gave
us relations~(iii), we have
$$ (uv^{-1})(u+v)^\oot=\q^{2/\tau}
   (u+v)^n(uv^{-1})=(u+v)^\oot(uv^{-1})\,.         $$
Thus, $(u+v)^\oot$ commutes with something that
is not (a series in) $u+v$, and therefore with
both~$u$ and~$v$ separately. Then it is a series
in $u^\oot$ and $v^\oot$, but the only such
series to scale right is $u^\oot+v^\oot$. Hence
\begin{equation*}\tag{iv}
   (u+v)^\oot=u^\oot+v^\oot,
\end{equation*}
and the rest is straightforward. By relation~(ii) and
the same argument as for the numerators we have
$$ (v^\oot+u^\oot;\q^{-2/\tau^2})_\infty
   =(v^\oot;\q^{-2/\tau^2})_\infty
   (u^\oot;\q^{-2/\tau^2})_\infty\,,               $$
then relation~(iv) turns this into
$$ ((u+v)^\oot;\q^{-2/\tau^2})_\infty
   =(v^\oot;\q^{-2/\tau^2})_\infty
   (u^\oot;\q^{-2/\tau^2})_\infty\,,               $$
and in their turn relations~(iii) allow to reunite
the numerators with denominators and obtain
$$ \x(u+v)=\x(u)\x(v)\,,                           $$
as we want. It only remains, therefore, to find out
if those hypothetical relations (ii-iv) actually hold
good if formal~$u$ and~$v$ are replaced back by
$$ u=\e^{-2\ipi z}\AND v=\e^{\tau\,d/dz}.          $$
But obviously\footnote
{See footnote to Lemma~2 below.}
$$ u^\oot=\e^{-2\ipi z/\tau}\AND v^\oot=\e^{d/dz}, $$
and therefore relation~(ii) holds as good as~(i),
relations~(iii) are checked trivially, and,
finally,~(iv) has already been established back in
Section~1 (relation~(\ref{oot})). So we are done
-- but another important point now needs
clearing up.

As we have just learnt, the numerators alone already
satisfy Sch\"utzen\-berger's equation -- and
with it in fact all the other noncommutative
identities in question and a full complement of so
called q-hypergeometric sums very similar to our
integral identities, only much older.\footnote
{See (Koornwinder 1996) for details and history.}
So, what we have actually shown so far is that
adding suitable denominators does no harm. But
what good does it do? The answer is already apparent
in Figure~1. Note that if $|\tau|=1$,\footnote
{Not to be mistaken for $|\q|=1$.}
then, on top of central symmetry, zeros and poles of
the $\g$-function are mirror symmetric to each other
about the line passing through~$1$ and~$\tau$, and
as a result $|\g(z)|=1$ everywhere on that line --
none of which can be said of the numerator because
it has no poles in the first place. So, if we are
to go beyond formal algebra and develop any kind of
a unitary theory, then, as it was first realized by
L. Faddeev, we really need the whole of the
$\g$-function, and $|\tau|=1$ is the case to look
into.\footnote
{In fact, the limit case when $\tau>0$ would do as
well, but we have to choose something.}
%%%%%%%%%%%%%%%%%%%%%%%%%%%%%%%%%%%%%%%%%%%%%%%%%%%%%
\section{\protect\large\bf
Case $\boldsymbol{|\tau|=1}$}
%%%%%%%%%%%%%%%%%%%%%%%%%%%%%%%%%%%%%%%%%%%%%%%%%%%%%
For convenience, let us adjust the `reference frame'
so that the aforementioned symmetry axis becomes the
real line. To this end, fix
``Planck's constant''~$\hbar$ and offset~$\opp$, set
$\om=\sqrt{-\pi\hbar/2\tau}$ and
$\op=\sqrt{-\pi\hbar\tau/2}$ --
so that conversely
$$ \hbar=-\,\frac{2\om\op}\pi\AND
   \tau=\frac\op\om\,,                             $$
and redefine the $\g$-function like this:
$$ \g^{}_{\mbox{new}}(z)=\g_{\smash{\mbox{old}}}
   \Bigl(\frac{z+\opp}{2\om}\Bigr)\,.              $$
In these terms, the original setup corresponds to
$\opp=0$ and $\hbar=-\,\tau/2\pi$, but now we opt
instead for $\opp=\om+\op$ and some positive~$\hbar$,
say, $\hbar=1/2\pi$ for a change.
%%%%%%%%%%%%%%%%%%%%%%%%%%%%%%%%%%%%%%%%%%%%%%%%%%%%%
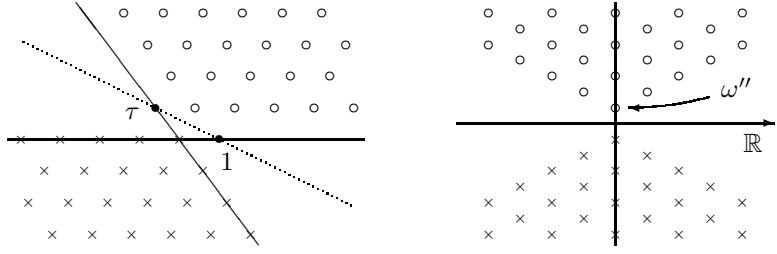
\begin{figure}[ht]\centering
%\framebox{
 \begin{picture}(300,92)
  \put(70,40){\begin{picture}(0,0)
   \put(30,-40){\line(-3,4){69}}
   \put(-65,0){\line(1,0){135}}
%   \put(70,-10){\makebox(0,0)[br]{$\IR$}}
   \put(-9,12){\circle*{3}}
   \put(-15,12){\makebox(0,0)[tr]{$\tau$}}
   \put(15,0){\circle*{3}}
   \put(18,-5){\makebox(0,0)[t]{$1$}}
   \dottedline{2}(65,-25)(-59,37)
%   \dottedline{2}(-20,-40)(26,52)
%   \put(65,-25){\line(-2,1){124}}
   \multiput(-21,48)(15,0){6}{\zero}
   \multiput(-12,36)(15,0){6}{\zero}
   \multiput(-3,24)(15,0){5}{\zero}
   \multiput(6,12)(15,0){5}{\zero}
   \multiput(0,0)(-15,0){5}{\pole}
   \multiput(9,-12)(-15,0){5}{\pole}
   \multiput(18,-24)(-15,0){6}{\pole}
   \multiput(27,-36)(-15,0){6}{\pole}
 \end{picture}}
  \put(235,46){\begin{picture}(0,0)
   \put(-60,0){\vector(1,0){120}}
   \put(0,-46){\line(0,1){92}}
   \put(52,-4){\makebox(0,0)[t]{$\IR$}}
   \put(5,6){\vector(-1,0){0}}
   \qbezier(5,6)(20,6)(35,10)
   \put(40,10){\makebox(0,0)[bl]{$\opp$}}
   \multiput(0,6)(0,12){4}{\zero}
   \multiput(12,12)(0,12){3}{\zero}
   \multiput(24,18)(0,12){3}{\zero}
   \multiput(36,24)(0,12){2}{\zero}
   \multiput(48,30)(0,12){2}{\zero}
   \multiput(-12,12)(0,12){3}{\zero}
   \multiput(-24,18)(0,12){3}{\zero}
   \multiput(-36,24)(0,12){2}{\zero}
   \multiput(-48,30)(0,12){2}{\zero}
   \multiput(0,-6)(0,-12){4}{\pole}
   \multiput(12,-12)(0,-12){3}{\pole}
   \multiput(24,-18)(0,-12){3}{\pole}
   \multiput(36,-24)(0,-12){2}{\pole}
   \multiput(48,-30)(0,-12){2}{\pole}
   \multiput(-12,-12)(0,-12){3}{\pole}
   \multiput(-24,-18)(0,-12){3}{\pole}
   \multiput(-36,-24)(0,-12){2}{\pole}
   \multiput(-48,-30)(0,-12){2}{\pole}
%   \dottedline{2}(55,6)(-55,-2)
 \end{picture}}\end{picture}
%}
\caption{Case~$|\tau|=1$ in old and new frames.}
\end{figure}
%%%%%%%%%%%%%%%%%%%%%%%%%%%%%%%%%%%%%%%%%%%%%%%%%%%%%
The infinite product expansion and defining equations
then read
\begin{gather}
   \g(z)=\frac{(1+\e^{-\ipi(z-\op)/\om})
   (1+\e^{-\ipi(z-3\op)/\om})
   (1+\e^{-\ipi(z-5\op)/\om})\ldots}
   {(1+\e^{-\ipi(z+\om)/\op})
   (1+\e^{-\ipi(z+3\om)/\op})
   (1+\e^{-\ipi(z+5\om)/\op})\ldots}\,,            \\
   \frac{\g(z+\op)}{\g(z-\op)}=1+\e^{-\ipi z/\om}
   \AND
   \frac{\g(z+\om)}{\g(z-\om)}=1+\e^{-\ipi z/\op},
\end{gather}
zeros/poles are located at the points $z=k\om+l\op$
with~$k$ and~$l$ positive/negative odd integers,
and $\g(z)\sim1$ as $z\rightarrow\infty$ in the
sector $|\arg z|<\arg\om=(\pi-\arg\tau)/2$ or in
particular as $z\rightarrow+\infty$. In their turn,
the reflection formula~(\ref{rf}), tau-gamma
integral~(\ref{ft}) and tau-binomial
theorem~(\ref{tbt}) become
\begin{align}
   \label{nrf}
   \g(z)\g(-z)&=\,\alpha\,\,\e^{\ipi z^2},        \\
   \g(z+\opp)&=\,\beta\rint\frac
   {\e^{-2\ipi\zeta z}}{\g(\opp-\zeta)}\,,        \\
   \label{ntbt}
   \frac{\g(y)\g(z+\opp)}{\g(y+z)}
   &=\,\beta\rint\frac{\e^{-2\ipi\zeta z}
   \g(y-\zeta)}{\g(\opp-\zeta)}\,,
\end{align}
where in both integral the real line is suitably
indented,\footnote
{See footnote to~(\ref{tbt}).}
and, if you must know,
$\alpha=-\,\tau\e^{-\ipi\opp{}^2}
\g(\om-\op)^2$ and $\beta=\g(\om-\op)/2\om$, and in
its turn, $\g(\om-\op)$ is what used to be~$\g(1)$
(see Sections~2 and~4). And, of course, on top of all
this we have the proto-unitarity property
$$ |\g(z)|=1\quad\text{for all $z\in\IR$,}         $$
which implies that operator $\g(A)$ is unitary
whenever~$A$ is self-adjoint, and allows to
restate our findings so far as follows.
Identities~(\ref{seven}) and~(\ref{pentagon}) become
\begin{theorem}
Let~$\Q$ and~$\P$ be (the self-adjoint closures in
$L^2(\IR)$ of) Schr\"o\-din\-ger's position and
momentum operators
$$ \Q f(z)=zf(z)\AND
   \P f(z)=\frac{f'(z)}{2\ipi}                     $$
or a unitary equivalent pair. Then operators
$\g(\pm\Q)$, $\g(\P)$ and $\g(\P+\Q)$ are all unitary
and satisfy the `$3=4$ identity'
\begin{equation}\label{nseven}
   \g(-\Q)\g(\P)\g(\Q)=\g(\P)\g(-\Q)\g(\Q)\g(\P)
\end{equation}
and `pentagon identity'
\begin{equation}\label{npentagon}
   \g(\P)\g(\Q)=\g(\Q)\g(\P+\Q)\g(\P)\,.
\end{equation}
\end{theorem}
Relation~(\ref{oot}) becomes
\begin{lemma}
Let~$u(t)$ and~$v(t)$ be Weyl--Stone--von--Neumann's
operators
$$ u(t)=\e^{2\ipi t\Q}\AND v(t)=\e^{2\ipi t\P},    $$
and let
$$ u\equiv u(2\op)=\e^{-\ipi\Q/\om}\AND
   v\equiv v(2\op)=\e^{-\ipi\P/\om}.               $$
Then
$$ u^\oot=u(2\om)=\e^{-\ipi\Q/\op}\AND
   v^\oot=v(2\om)=\e^{-\ipi\P/\op},                $$
and furthermore
\begin{equation}
   (u+v)^\oot=u^\oot+v^\oot,
\end{equation}
provided the branch is so chosen that
$(\e^{-\ipi z/\om})^\oot=\e^{-\ipi z/\op}$ for all
$z\in\IR$.\footnote
{Fittingly, such a branch only fails to exist
if~$\om$ is real -- which it is absolutely not.}
\end{lemma}
Finally, identities~(\ref{sb}) and~(\ref{cc}) become
\begin{theorem}
If function~$\x$ is such that
$\x(\e^{-\ipi z/\om})=\g(z)$ for all $z\in\IR$,
that is
$$ \x(w)=\frac{(1+\q w)(1+\q^3w)(1+\q^5w)\ldots}
   {\bigl(1+\q^{-1/\tau^2}w^\oot\bigr)
   \bigl(1+\q^{-3/\tau^2}w^\oot\bigr)
   \bigl(1+\q^{-5/\tau^2}w^\oot\bigr)\ldots}       $$
with the same proviso about the branch, then
operators $\x(u)$, $\x(v)$, $\x(u+v)$ and
$\x(u+\q vu+v)$ are all
unitary and satisfy `Sch\"utzenberger's identity'
\begin{equation}\label{nsb}
   \x(u+v)=\x(u)\x(v)
\end{equation}
and the other way around identity
\begin{equation}
   \x(v)\x(u)=\x(u+\q vu+v)\,.
\end{equation}
\end{theorem}
It is straightforward to upgrade the arguments of
Section~6 to the level of strict proofs.\footnote
{See also (Woronowicz 2000) and
(Bytsko Teschner 2003) for alternative takes
on the subject.}
It should be noted though that it was mostly for
demonstration purposes that in that Section every
identity was independently derived all the way from
the tau-binomial theorem. It would be more practical
to derive only one of them and then transform it
into the remaining three by purely `noncommutative'
techniques. For instance, the pentagon
identity~(\ref{npentagon}) can be easily transformed
into $3=4$ by either of the following ways.

Recall
\begin{lemma}[folklore]
Operators
$$ \sigma_1=\alpha\,\e^{i\Q^2/2\hbar}\AND
   \sigma_2=\alpha\,\e^{i\P^2/2\hbar}              $$
satisfy Artin's braid group relation\footnote
{By the way, these two triple products equal not only
each other in fact, but also ($\sqrt i$ times)
Fourier transform understood as a unitary operator
in~$L^2(\IR)$, that is
$\sqrt i\,\e^{\ipi (\Q^2+\P^2-1/2\pi)}$. I leave it
to the reader to figure this out.}
\begin{equation}\label{artin}
   \sigma_1\sigma_2\sigma_1
   =\sigma_2\sigma_1\sigma_2\,.
\end{equation}
\end{lemma}
\begin{proof}
By the product differentiation rule we have
$$ \P\,\e^{i\Q^2/2\hbar}
   =\e^{i\Q^2/2\hbar}\,(\P+\Q)\,,                  $$
and by unitary equivalence
$$ (\P+\Q)\,\e^{i\P^2/2\hbar}
   =\e^{i\P^2/2\hbar}\,\Q\,.                       $$
Hence
$$ \P\,\e^{i\Q^2/2\hbar}\e^{i\P^2/2\hbar}
   =\e^{i\Q^2/2\hbar}\e^{i\P^2/2\hbar}\,\Q\,,      $$
and the result follows at once.
\end{proof}
Now, `divide' Artin's relation by the pentagon
identity~(\ref{npentagon}),\footnote
{I write $A/B$ for $AB^{-1}$ whenever~$A$ and~$B$
commute.}
$$ \frac1{\g(\Q)}\,\frac1{\g(\P)}\,\sigma_1\sigma_2
   \sigma_1=\frac1{\g(\P)}\,\frac1{\g(\P+\Q)}\,
   \frac1{\g(\Q)}\,\sigma_2\sigma_1\sigma_2\,,     $$
then use the product differentiation rule (see the
above proof) and reflection formula~(\ref{nrf}) to
simplify the left hand side like this:
\begin{multline*}
   \alpha^3\,\frac1{\g(\Q)}\,\frac1{\g(\P)}\,
   \e^{i\Q^2/2\hbar}\e^{i\P^2/2\hbar}
   \e^{i\Q^2/2\hbar}                               \\
   =\alpha^3\,\frac{\e^{i\Q^2/2\hbar}}{\g(\Q)}\,
   \e^{i\P^2/2\hbar}\,
   \frac{\e^{i\Q^2/2\hbar}}{\g(\Q)}
   =\g(-\Q)\g(\P)\g(-\P)\g(-\Q)\,,
\end{multline*}
and the right hand side like this:
\begin{multline*}
   \alpha^3\,\frac1{\g(\P)}\,\frac1{\g(\P+\Q)}\,
   \frac1{\g(\Q)}\,\e^{i\P^2/2\hbar}
   \e^{i\Q^2/2\hbar}\e^{i\P^2/2\hbar}              \\
   =\alpha^3\,\frac{\e^{i\P^2/2\hbar}}{\g(\P)}\,
   \frac{\e^{i\Q^2/2\hbar}}{\g(\Q)}\,
   \frac{\e^{i\P^2/2\hbar}}{\g(-\P)}
   =\g(-\P)\g(-\Q)\g(\P)\,.
\end{multline*}
Hence
$\g(-\P)\g(-\Q)\g(\P)=\g(-\Q)\g(\P)\g(-\P)\g(-\Q)$,
which is the $3=4$ identity~(\ref{seven}) modulo
unitary equivalence. So, one way to transform the
pentagon identity~(\ref{npentagon}) into $3=4$, or
vise versa for that matter, is use the formula
$$ (\text{$3=4$ identity})
   =(\text{pentagon identity})^{-1}
   (\text{Artin's relation})\,.                    $$

The other transformation is, in contrast, one way
only, and it goes like this:
\begin{multline*}
   \g(-\Q)\underline{\g(\P)\g(\Q)}
   =\g(-\Q)\g(\Q)\g(\P+\Q)\g(\P)
   =\g(\Q)\underline{\g(-\Q)\g(\P+\Q)}\g(\P)       \\
   =\underline{\g(\Q)\g(\P+\Q)\g(\P)}\g(-\Q)
   =\g(\P)\g(-\Q)\g(\Q)\g(\P)\,.
\end{multline*}
This time all is done with the pentagon relation,
which is first used `as is', then in its unitary
equivalent form
$$ \g(-\Q)\g(\P+\Q)=\g(\P+\Q)\g(\P)\g(-\Q),        $$
and then again `as is'.
%%%%%%%%%%%%%%%%%%%%%%%%%%%%%%%%%%%%%%%%%%%%%%%%%%%%%
\section{\protect\large\bf Yang-Baxterization}
%%%%%%%%%%%%%%%%%%%%%%%%%%%%%%%%%%%%%%%%%%%%%%%%%%%%%
The bottom line so far is that we have obtained
every operator interpretation of the tau-binomial
theorem~(\ref{tbt}) I know of. It only remains, then,
to do the same to the ultimate integral
identity~(\ref{pfaff}). Unfortunately, due to the
greater number of variables involved, there are
more such interpretations than would be appropriate
in an introductory article. We will, therefore, leave
Theorem~3 for another time and limit ourselves to
generalization of Theorem~1 and Lemma~4. Here it is. 
\begin{theorem}
For all $\lambda,\mu\in\IR$ there hold
quasi-Yang-Baxter equations
\begin{multline}\label{ybp}
   \g(\P)\g(\lambda-\P)\g(\mu+\P-\Q)
   \g(\lambda-\P+\Q)\g(\Q)\g(\mu-\Q)               \\
   =\g(\mu-\Q)\g(\Q)\g(\P+\Q)
   \g(\lambda+\mu-\P-\Q)\g(\P)\g(\lambda-\P)\,
\end{multline}
\begin{equation}\label{ybs}
   \frac{\g(-\Q)}{\g(\lambda-\Q)}\,\frac{\g(\P)}
   {\g(\lambda+\mu+\P)}\,\frac{\g(\Q)}{\g(\mu+\Q)}
   =\frac{\g(\P)}{\g(\mu+\P)}\,\frac{\g(-\Q)\g(\Q)}
   {\g(\lambda-\Q)\g(\mu+\Q)}\,
   \frac{\g(\P)}{\g(\lambda+\P)}
\end{equation}
and the true Yang-Baxter equation
\begin{equation}
   \sigma_1(\lambda)\,\sigma_2(\lambda+\mu)\,
   \sigma_1(\mu)=\sigma_2(\mu)\,
   \sigma_1(\lambda+\mu)\,\sigma_2(\lambda)\,,
\end{equation}
where~$\sigma(\lambda)$ is Fateev-Zamolodchikov's
R-matrix:\footnote
{It is called that after its finite-dimensional
relative from (Fateev Zamolodchikov 1982).}
$$ \sigma_1(\lambda)=\frac{\sigma_1}
   {\g(\frac\lambda2+\Q)\g(\frac\lambda2-\Q)}\AND
   \sigma_2(\lambda)=\frac{\sigma_2}
   {\g(\frac\lambda2+\P)\g(\frac\lambda2-\P)}\,.   $$
\end{theorem}
\begin{proof}
With some patience, all three identities could be
derived starting from the ultimate integral
identity~(\ref{pfaff}) and following the guidelines
of Section~6, which is left as another exercise.
This would not quite prove the theorem though, for,
as we remember, the said integral identity has not
been actually verified. We need, therefore, some
kind of a direct `noncommutative' proof, and this is
where the techniques shown in the previous section
come into their own.

Apply the pentagon identity to the underlined
pieces either as is or in a suitable unitary
equivalent form:
\begin{multline*}
   \g(\P)\underline{\g(\lambda-\P)\g(\mu+\P-\Q)}
   \g(\lambda-\P+\Q)\g(\Q)\g(\mu-\Q)               \\
   =\g(\P)\g(\mu+\P-\Q)\g(\lambda+\mu-\Q)
   \underline{\g(\lambda-\P)
   \g(\lambda-\P+\Q)\g(\Q)}\g(\mu-\Q)              \\
   =\g(\P)\g(\mu+\P-\Q)\g(\lambda+\mu-\Q)\g(\Q)
   \underline{\g(\lambda-\P)\g(\mu-\Q)}            \\
   =\underline{\g(\P)\g(\mu+\P-\Q)\g(\mu-\Q)}
   \g(\Q)\g(\lambda+\mu-\Q)
   \g(\lambda+\mu-\P-\Q)\g(\lambda-\P)             \\
   =\g(\mu-\Q)\underline{\g(\P)\g(\Q)}\g(\lambda
   +\mu-\Q)\g(\lambda+\mu-\P-\Q)\g(\lambda-\P)     \\
   =\g(\mu-\Q)\g(\Q)\g(\P+\Q)
   \underline{\g(\P)\g(\lambda+\mu-\Q)
   \g(\lambda+\mu-\P-\Q)}\g(\lambda-\P)            \\
   =\g(\mu-\Q)\g(\Q)\g(\P+\Q)\g(\lambda+\mu-\P-\Q)
   \g(\P)\g(\lambda-\P)\,.
\end{multline*}
This settles~(\ref{ybp}), and then~(\ref{ybs})
follows in exactly the same way as in the last
section the $3=4$ followed from the pentagon
identity, that is by `dividing' Artin's
relation~(\ref{artin}) by~(\ref{ybp}).
In its turn, the Yang-Baxter equation emerges
if~(\ref{artin}) is divided, instead of~(\ref{ybp}),
by its unitary equivalent variant
\begin{multline*}
   \g(\tfrac\lambda2+\P)\g(\tfrac\lambda2-\P)
   \g(\tfrac{\lambda+\mu}2+\P-\Q)
   \g(\tfrac{\lambda+\mu}2-\P+\Q)
   \g(\tfrac\mu2+\Q)\g(\tfrac\mu2-\Q)              \\
   =\g(\tfrac\mu2-\Q)\g(\tfrac\mu2+\Q)
   \g(\tfrac{\lambda+\mu}2+\P+\Q)
   \g(\tfrac{\lambda+\mu}2-\P-\Q)
   \g(\tfrac\lambda2+\P)\g(\tfrac\lambda2-\P)\,.
\end{multline*}
And finally, if you have done the exercise suggested
at the beginning of the proof, you can reverse it and
thus settle~(\ref{pfaff}).
\end{proof}

To conclude, I want to thank
R.~Kashaev,
I.~Loris,
V.~Matveev,
Yu.~Melnikov,
M.~Se\-me\-nov-Tian-Shan\-sky
and S.~Shkarin
for helpful discussion,
A.~Alekseev and F.~Lambert for support,
and L.~Faddeev for patience.
%%%%%%%%%%%%%%%%%%%%%%%%%%%%%%%%%%%%%%%%%%%%%%%%%%%%%
\section*{\protect\large\bf References}
\noindent E.W. Barnes.
\em The genesis of the double gamma function,
\em Proc. London Math. Soc. {\bf 31} (1899) 358-381
\\[1ex]A.G. Bytsko and J. Teschner.
\em R-operator, co-product and Haar-measure for the
modular double of $U_q(\mathfrak{sl}(2,\mathsf{R})$,
\em Comm. Math. Phys. {\bf 240} (2003) 171-196\\
{}[math.QA/0208191]
\\[1ex]L. Faddeev.
\em Discrete Heisenberg-Weyl group and modular group,
\em Lett. Math. Phys. {\bf 34} (1995) 249-254
[hep-th/9504111]
\\[1ex]L. Faddeev.
\em Modular double of a quantum group,
\em Math. Phys. Stud. {\bf 21} (2000) 149-156
[math.QA/9912078]
\\[1ex]L. Faddeev, R. Kashaev and
A.Yu. Volkov.
\em Strongly coupled quantum discrete Liouville
theory. I: Algebraic approach and duality,
\em Comm. Math. Phys. {\bf 219} (2001) 199-219
[hep-th/0006156]
\\[1ex]V. Fateev and A. Zamolodchikov.
\em Selfdual solutions of the star triangle relations
in $Z(N)$ models,
\em Phys. Lett. {\bf A92} (1982) 37-39
\\[1ex]R. Kashaev,
\em On the spectrum of Dehn twists in quantum
Teichmuller theory\\
\em[math.QA/0008148]
\\[1ex]S. Kharchev, D. Lebedev and
M. Semenov-Tian-Shansky.
\em Unitary representations of
$U_q(\mathfrak{sl}(2,\mathsf{R})$, the modular double,
and the multiparticle q-deformed Toda chains,
\em Comm. Math. Phys. {\bf 225} (2003) 573-609
[hep-th/0102180]
\\[1ex]T. Koornwinder.
\em Special functions and q-commuting variables.
\em[q-alg/9608008]
\\[1ex]B. Ponsot and J. Teschner.
\em Clebsch-Gordan and Racah-Wigner coefficients for
a continuous series of representations of
$U_q(\mathfrak{sl}(2,\mathsf{R})$,
\em Comm. Math. Phys. {\bf 224} (2001) 613-655
[math.QA/0007097]
\\[1ex]M.-P. Sch\"utzenberger.
\em Une interpr\`etation de certaines solutions de
l'\`equa\-tion fonctionelle: $F(x+y)=F(x)F(y)$,
\em C. R. Acad. Sci. Paris {\bf 236} (1953) 352-353
\\[1ex]T. Shintani.
\em On a Kronecker limit formula for real quadratic
fields,
\em J. Fac. Sci. Univ. Tokyo Sect. 1A Math. {\bf 24}
(1977) 167-199
\\[1ex]A. Yu. Volkov.
\em Beyond the `Pentagon identity'
\em Lett. Math. Phys. {\bf 39} (1997) 393-397
[q-alg/9603003]
\\[1ex]S.L. Woronowicz.
\em Quantum exponential function,
\em Rev. Math. Phys. {\bf 136} (2000) 873-920
%%%%%%%%%%%%%%%%%%%%%%%%%%%%%%%%%%%%%%%%%%%%%%%%%%%%%
\end{document}